\documentclass[a4paper,reqno]{amsart}
\usepackage{amssymb}
\usepackage{latexsym}
\usepackage{amsmath}
\usepackage{euscript}
\usepackage{graphics,color}
\usepackage[all]{xy}
\usepackage[margin=3cm]{geometry}
\usepackage{MnSymbol} 

      \def\dR{{\mathbb R}}

   \def\dZ{{\mathbb Z}}

\newcommand{\be}{\begin{equation}}
\newcommand{\ee}{\end{equation}}
\newcommand{\ba}{\begin{eqnarray}}
\newcommand{\ea}{\end{eqnarray}}
\newcommand{\baa}{\begin{eqnarray*}}
\newcommand{\eaa}{\end{eqnarray*}}
\newcommand{\bb}{\begin{thebibliography}}
\newcommand{\eb}{\end{thebibliography}}

\newcommand{\bi}[1]{\bibitem{#1}}
\newcommand{\lab}[1]{\label{#1}}
\newcommand{\re}[1]{(\ref{#1})}

\newcounter{my}
\newcommand{\he}%
   {\stepcounter{equation}\setcounter{my}%
   {\value{equation}}\setcounter{equation}0%
   }%
\newcommand{\she}%
   {\setcounter{equation}{\value{my}}%
    }%

\renewcommand\t{\tilde}

\newtheorem{theorem}{Theorem}[section]

\theoremstyle{definition}

\newtheorem{remark}[theorem]{Remark}

\numberwithin{equation}{section}

\begin{document}

\title{An algebraic description of the bispectrality of the biorthogonal rational functions of Hahn type}

\author{Satoshi Tsujimoto}
\author{Luc Vinet}
\author{Alexei Zhedanov}

\address{Department of Applied Mathematics and Physics, 
Graduate School of Informatics, Kyoto University,
Yoshida-Honmachi, Kyoto, Japan 606-8501
}

\address{Centre de recherches math\'ematiques,
Universit\'e de Montr\'eal, P.O. Box 6128, Centre-ville Station,
Montr\'eal (Qu\'ebec), Canada H3C 3J7}

\address{School of Mathematics, Renmin University of China, Beijing 100872, China}

\begin{abstract}
The biorthogonal rational functions of the ${_3}F_2$ type on the uniform grid provide the simplest example of rational functions with bispectrality properties that are similar to those of classical orthogonal polynomials. These properties are described by three difference operators $X,Y,Z$ which are tridiagonal with respect to three distinct bases of the relevant finite-dimensional space.  The pairwise commutators of the operators $X,Y,Z$ generate a quadratic algebra which is akin to the algebras of Askey-Wilson type attached to hypergeometric polynomials.
\end{abstract}



\maketitle

\section{Introduction}

Among the families of biorthogonal rational functions \cite{W}, \cite{GM}, \cite{IM}, some classes of particular interest are those that are bispectral namely, those where the functions satisfy both a recurrence relation and a differential/difference equation. This paper explores the algebraic characterization of this feature for biorthogonal rational functions through the study of those of Hahn type \cite{Zhe_RIMS}. 
We may remark that related $q$-rational functions have been studied in \cite{KM} using symmetry techniques.

The algebraic description of bispectrality in the case of orthogonal polynomials is well developed especially for the families that belong to the Askey scheme \cite{KLS} and its Bannai-Ito extension \cite{BI}. It centers around the Askey-Wilson algebra \cite{GLZ_Annals} and its contractions and specializations. These algebras are realized by the bispectral operators in either the variable or the degree representation. In finite dimensions these operators form the so-called Leonard pairs which underpins the Leonard duality of the polynomials of the Askey scheme. Let us recall that a Leonard pair \cite{Ter} roughly amounts to two linear operators that are each tridiagonal in an eigenspace of the other and that
polynomials enjoying the Leonard duality \cite{L} are such that their duals, obtained by exchanging the variable and the degree, are also orthogonal.

Biorthogonal functions naturally arise from generalized eigenvalue problems (GEVPs) of the form $L_1U(x)=\lambda L_2 U(x)$ where $L_1$ and $L_2$ are two operators acting on functions of a variable $x$ and $\lambda$ is the eigenvalue. It is understood that rational functions occur as solutions of the GEVP \cite{Zhe_GEVP} when the two operators $L_1$ and $L_2$ act tridiagonally in a certain basis. We wish to develop the understanding of bispectral situations in this context. From the study of a simple case, we contend that bispectrality rests on the presence of three operators to be denoted $X$, $Y$ and $Z$ that act tridiagonally on the given set of biorthogonal rational functions. In the case at hand, we shall observe that the difference equation and the recurrence relation emerge from the GEVPs associated respectively to the pairs $(X, Y)$ and $(X,Z)$. This triplet, albeit with distinct features, is hence playing for bispectral biorthogonal rational functions a role analogous to the one of Leonard pairs for orthogonal polynomials. The algebra realized by the triplet $X$, $Y$, $Z$ attached to the ${_3}F_2$ rational functions is quadratic and poised to encode the bispectral properties of these biorthogonal set much like the Askey-Wilson algebra and its specializations do for the hypergeometric orthogonal polynomials.

The paper will unfold as follows. The triplet of operators $X$, $Y$, $Z$ is introduced next in Section 2. It will be explained that these belong to the class of rational Heun operators on the uniform lattice. The rational functions $U_n(x;\alpha,\beta,N)$ of the ${_3}F_2$ type will be explicitly obtained in Section 3 as solutions of the GEVP involving $X$ and $Y$. This will be achieved using the action of these operators in a basis defined by ratios of two Pochhammer symbols. That $X$, $Y$ and $Z$ all generate shifts by one of the parameter $\alpha$ in $U_n(x;\alpha,\beta,N)$ will then be observed. The biorthogonal partners $V_n(x;\alpha,\beta,N)$ will be identified in Section 4 with the help of the symmetry properties of the hypergeometric weight distribution. Section 5 will provide the tridiagonal action of the triplet $(X, Y, Z)$ on the set of functions $U_n(x;\alpha,\beta,N)$. How their difference equation and recurrence relation are obtained from GEVPs associated to the triplet $(X, Y, Z)$ is the object of Section 6. The algebra realized by these operators is described in Section 7 and Section 8 will consist of concluding remarks.

\section{A triplet of difference operators}
Let us introduce the difference operators
\ba & X^{(\alpha , \beta)}= (x-\alpha) \mathcal{I} - x T^{-}, \lab{X_op}
\\
& Z^{(\alpha , \beta)}= -\mathcal{I} + \dfrac{x}{x-\alpha} T^{-} = (\alpha-x)^{-1} X,  \lab{Z_op} 
\\
& Y^{(\alpha , \beta)}=A_1(x)T^{+} + A_2(x) T^{-}+ A_0(x) \mathcal{I} \lab{Y_op} \ea
where
\ba
&&A_1(x) = (x-\alpha)(x-N)(x+ 1-\alpha), \; A_2(x) = x(x-\alpha) (x+\beta-\alpha-N), \nonumber \\
&&A_0(x) = (x-\alpha) \left( -2\,{x}^{2}+ \left( 2\,\alpha-1+2\,N-\beta \right) x-N \left( \alpha-1
 \right) \right). \lab{AY_op} \ea
 Note that $Y^{(\alpha , \beta)}$ has $(x-\alpha)$ as a global factor.
In \re{X_op}-\re{Y_op} we use the notation 
\be
T^{\pm} f(x) = f(x \pm 1), \quad \mathcal{I} f(x) = f(x). \lab{TI_def} \ee
These operators act on the space $M$ of real functions $f(x)$ defined on the linear grid $x=0, \pm 1, \pm 2, \dots$; that is $M = \{f: \dZ \rightarrow \dR\}$.
In what follows we will restrict $M$ to the space $M_N$ of functions defined on the finite set of points $x=0,1,\dots,N$.  $M_N$ is  hence of dimension $N+1$ and one has the following standard basis on this space:
\be
e_k(x) = \delta_{kx}, \quad k,x =0,1,2, \dots, N \lab{e_basis} \ee
In this basis the operators $X,Y,Z$ become matrices of size $(N+1) \times (N+1)$. More precisely, the matrices $X$ and $Z$ are lower bi-diagonal
\be
X=\left[
\begin{array}{ccccc}
-\alpha & 0 & 0 & \cdots & 0 \\
-1 & 1-\alpha & 0 &  & \vdots  \\
0 & \ddots & \ddots & \ddots & 0 \\
\vdots & & 1-N & N-\alpha-1 & 0 \\
0 & \cdots & 0 & -N & N-\alpha
\end{array}
\right],\lab{X_e} \ee

\be
Z=\left[
\begin{array}{ccccc}
-1 & 0 & 0 & \cdots & 0 \\
(1-\alpha)^{-1} & -1 & 0 &  & \vdots  \\
0 & \ddots & \ddots & \ddots & 0 \\
\vdots & & (N-1)(N-1-\alpha)^{-1} & -1 & 0 \\
0 & \cdots & 0 & N(N-\alpha)^{-1} & -1
\end{array}
\right]. 
\lab{Z_e} \ee
while the matrix $Y$ is tridiagonal
\be
Y=\left[
\begin{array}{ccccc}
A_0(0) & A_1(0) & 0 & \cdots & 0 \\
A_2(1) & A_0(1) & A_1(1) &  & \vdots  \\
0 & \ddots & \ddots & \ddots & 0 \\
\vdots & & A_2(N-1) & A_0(N-1) & A_1(N-1) \\
0 & \cdots & 0 & A_2(N) & A_0(N)
\end{array}
\right].
\lab{Y_e} \ee

As indicated in the Introduction, the three operators $X, Y, Z$ will be central in our study. In the spirit of our recent analysis of rational Heun operators \cite{TVZ}, they can be viewed as belonging to the set of operators that satisfy the rational Heun property on the linear grid. Let us elaborate a little on this remark.

Consider the rational functions $R_n(x), \: n=0,1,\dots $ of type $[n/n]$ in $M_N$ with poles at the prescribed locations $a_1, a_2, \dots, a_n$ on the real axis:
\be
R_0(x) = 1, \quad R_1(x) = r_{10} +r_{11} (x-a_1)^{-1}, \dots  R_n(x) = \xi_{n0} + \sum_{s=1}^n r_{ns} (x-a_s)^{-1}  
\lab{space_R} \ee
where $r_{ns}$ are some coefficients.

The abbreviation $[n/n]$ means that any such function can be presented as a ratio of two polynomials of degree $n$
\be
R_n(x) = \frac{Q_n(x)}{\Omega_n(x)}, \lab{RQP} \ee 
where the polynomial $\Omega_n(x)$ is the characteristic polynomial of the set $\{a_1,a_2, \dots, a_n\}$. The polynomials $Q_n(x)$ and $\Omega_n(x)$ need not be coprime. Indeed, if some coefficients $r_{ns}$ vanish in \re{space_R} then the rational function $R_n(x)$ has a smaller number of poles than $n$. And, this leads to the conclusion that the polynomials $Q_n(x)$ and $\Omega_n(x)$ have co coinciding zeros (on $\dR$).  

The rational Heun operators were introduced in \cite{TVZ} as operators $W$ transforming any rational function of type $[n/n]$ with poles at prescribed locations $a_1, a_2, \dots, a_n$ to a function of type $[(n+1)/(n+1)]$ with poles at $a_1, a_2, \dots, a_n, a_{n+1}$. Roughly speaking these operators exhibit a ``raising" property according to which a pole at $a_{n+1}$ is added to those at $a_1,a_2, \dots, a_n$.

The rational Heun operators acting on the Askey-Wilson grid were constructed in \cite{TVZ} (see also \cite{VZ}) and seen to lead to the Wilson biorthogonal functions of type ${_{10}}\Phi_9$ \cite{W}. Here we consider the simplest case of the linear grid with poles located at $ \alpha+n-1, \: n=1,2, \dots$ .The second-order operators in that class are of the form
\be
W=A_1(x) T^+ + A_2(x) T^- + A_0(x) \mathcal{I}, \lab{W_gen} 
\ee
with the functions $A_0(x),A_1(x),A_2(x)$ such that $W$ transforms any elementary rational monomial $(x-\alpha-n)^{-1}$ into a linear combination of the monomials $1, (x-\alpha)^{-1}, (x-\alpha-1)^{-1}, \dots, (x-\alpha-n-1)^{-1}$. When $A_1(x) =0$, $W$ is in fact a first-order difference operator. 
It is easily verified that the operators $X^{(\alpha,\beta)}$, $Z^{(\alpha,\beta)}$ and $Y^{(\alpha,\beta)}$ given in \re{X_op}, \re{Y_op} and \re{Z_op} respectively, all meet this criterion.
More precisely we have for $n=0, 1,2,\dots$
\ba 
& X^{(\alpha,\beta)} \{(x-\alpha-n)^{-1}\} = {\dfrac {n}{x-\alpha-n}}- {\dfrac {1+\alpha+n}{x-1-\alpha-n}}, \lab{X_mon} 
\\
& Z^{(\alpha,\beta)} \{(x-\alpha-n)^{-1}\} = -{\dfrac {\alpha}{ \left( 1+n \right)  \left( x-\alpha \right) }}-
 \dfrac{1}{x-\alpha-n} +{\dfrac {1+\alpha+n}{ \left( 1+n
 \right)  \left( x-1-\alpha-n \right) }} \lab{Z_mon} \ea
and
\be
Y^{(\alpha,\beta)} \{(x-\alpha-n)^{-1}\} = \beta+1 + \frac{\kappa_n^{(1)}}{x-\alpha-n+1} + \frac{\kappa_n^{(2)}}{x-\alpha-n} + \frac{\kappa_n^{(3)}}{x-\alpha-n-1}, \lab{Y_mon} \ee
where the residues $\kappa_n^{(1)},\kappa_n^{(2)},\kappa_n^{(3)}$ are cubic polynomials in $n$.

As shall be seen these special rational Heun operators possess remarkable factorization properties in bases of interest. We will immediately call upon them to set up the GEVP that will possess as solutions the biorthogonal rational functions we wish to concentrate on.

\section{The rational functions of the ${_3}F_2$ type as solutions of a generalized eigenvalue problem} \lab{secrat}

We wish to solve the GEVP:
\be
Y^{(\alpha , \beta)}U=\lambda X^{(\alpha , \beta)}U.\lab{GEP}
\ee
To that end we introduce the following basis, again on the space $M_N$ of functions defined on the finite set of points $x=0,1,\dots, N$:
\be
\varphi_n(x;\alpha) = \frac{(-x)_n}{(\alpha-x)_n}, \quad n=0,1,\dots, N ,  \lab{phi_n} \ee where
$(a)_n = a(a+1) \dots (a+n-1)$ is the Pochhammer symbol (or the shifted factorial). Enforcing the restriction to $M_N$,
it is verified that the operators $X^{(\alpha , \beta)},Y^{(\alpha , \beta)},Z^{(\alpha , \beta)}$ act in a tridiagonal fashion in this basis:
\ba
&&X^{(\alpha , \beta)} \varphi_n(x;\alpha) = -n  \varphi_{n+1}(x;\alpha) + (n-\alpha) \varphi_n(x;\alpha) ,\quad 0\le n<N,  \\ \nonumber
&& X^{(\alpha , \beta)} \varphi_N(x;\alpha) =  (N-\alpha) \varphi_N(x;\alpha) \lab{X_phi} \ea 

\noindent

\ba
&&Z^{(\alpha , \beta)} \varphi_n(x;\alpha) = \varphi_{n+1}(x;\alpha) - \varphi_n(x;\alpha), \quad 0\le n<N,  \\    \nonumber
&&Z^{(\alpha , \beta)} \varphi_N(x;\alpha) =  - \varphi_N(x;\alpha)
\ea  \lab{Z_phi} 
and
\ba 
&&Y^{(\alpha , \beta)} \varphi_n(x;\alpha) = \nu_n^{(1)} \varphi_{n+1}(x;\alpha)  + \nu_n^{(2)}\varphi_n(x;\alpha) + \nu_n^{(3)}\varphi_{n-1}(x;\alpha),    \quad 0\le n<N   \\ \nonumber 
&&Y^{(\alpha , \beta)} \varphi_N(x;\alpha) =  \nu_N^{(2)}\varphi_N(x;\alpha) + \nu_N^{(3)}\varphi_{N-1}(x;\alpha)
\lab{Y_phi} \ea
where
\ba
&&\nu_n^{(1)}= n^2 (n+\beta-N), \nonumber\\ 
&&\nu_n^{(2)} = -2n^3 +[2(1+N) +\alpha-\beta]n^2 +[-N-1+\alpha(\beta-N)]n, \nonumber \\
&&
\nu_n^{(3)}=n(N-n+1)(\alpha-n+1) \lab{nu_3}. \ea 
Equivalently, these formulas mean that the operators $X,Y,Z$ are given in the basis $\varphi_n(x)$ by the matrices 
\be
X=\left[
\begin{array}{ccccccc}
-\alpha & 0 & 0 & \cdots & 0 \\
0 & 1-\alpha & 0 &  & \vdots  \\
0 & -1 & 2-\alpha & \ddots &  0 \\
 & \ddots & \ddots & \ddots & 0 \\
\vdots & & 2-N & N-\alpha-1 & 0 \\
0 & \cdots & 0 & 1-N & N-\alpha
\end{array}
\right],\lab{XM_phi} \ee

\be
Y=\left[
\begin{array}{ccccc}
\nu_0^{(2)} & \nu_1^{(3)} & 0 & \cdots & 0 \\
\nu_0^{(1)} & \nu_1^{(2)} & \nu_2^{(3)} &  & \vdots  \\
0 & \ddots & \ddots & \ddots & 0 \\
\vdots & & \nu_{N-2}^{(1)} & \nu_{N-1}^{(2)} & \nu_N^{(3)} \\
0 & \cdots & 0 & \nu_{N-1}^{(1)} & \nu_N^{(2)}
\end{array}
\right],
\lab{YM_phi} \ee

\be
Z=\left[
\begin{array}{ccccc}
-1 & 0 & 0 & \cdots & 0 \\
1 & -1 & 0 &  & \vdots  \\
0 & \ddots & \ddots & \ddots & 0 \\
\vdots & & 1 & -1 & 0 \\
0 & \cdots & 0 & 1 & -1
\end{array}
\right]. 
\lab{ZM_phi} \ee

Note that in the basis formed by the functions $\varphi_n(x;\alpha)$, the operator $Y^{(\alpha , \beta)}$ admits the factorization
\be
Y^{(\alpha , \beta)} = X^{(\alpha , \beta)}\mathcal{Y}, \lab{Y-Y1} \ee 
where the matrix $\mathcal{Y}$ is bi-diagonal
\be
\mathcal{Y} \varphi_n(x;\alpha) = n(N-\beta-n)\varphi_n(x;\alpha) + n(n-N-1)\varphi_{n-1}(x;\alpha). \lab{Y_1} \ee
This property allows for an explicit hypergeometric solution of the GEVP \re{GEP}. Indeed, assume that  $U_n(x)$ is expanded
over the set of functions $\{\varphi_k(x;\alpha), k=0,\dots,n\}$.  In that basis, thanks to \re{Y-Y1}, we can factor out $X^{(\alpha , \beta)}$ on both sides of \re{GEP} to convert the GEVP into the standard eigenvalue equation
\be
\mathcal{Y}U_n=\lambda_nU_n\lab{GEVP}
\ee
with  the eigenvalue $\lambda _n$ is given by
\be
\lambda _n= n(N-\beta-n)\lab{EigV}
\ee
in view of \re{Y_1}. Now write
\be
U_n=\sum_{k=0}^n C_{n,k} \varphi_k(x;\alpha) \lab{U_phi},
\ee
it is readily seen using the action \re{Y_1} of $\mathcal{Y}$ on the basis functions $\varphi_k(x;\alpha)$ that \re{GEVP} implies the following
two-term recurrence relation for the coefficients $C_{n,k}$:
\be
(k+1)(k-N)C_{n,k+1}=\left[n(N-\beta -n)-k(N-\beta-k)\right]C_{n,k}
\ee
whose solution is found to be
\be
C_{n,k}=C_{n,0}\frac{(-n)_k(\beta+n-N)_k}{k!(-N)_k}.
\ee
Substituting into \re{U_phi} and using the definition \re{phi_n} of the functions $\varphi_k(x;\alpha)$, we thus arrive at the following solutions of the GEVP \re{GEP}\footnote{Note that when $n=N$ in \eqref{HIF} the sum defining the series runs up to $N$ in spite of the cancellation so as to yield as should be a rational function of type $[N/N]$.}:
\be
U_n(x;\alpha,\beta,N)= \frac{(-1)^n(-N)_n}{(\beta+1)_n} \: {_3}F_2\left({-x,-n,\beta+n-N\atop -N,\alpha-x}; 1 \right ) 
\lab{HIF} 
\ee
with eigenvalues $\lambda_n$ given by \re{EigV}. The coefficient $C_{n,0}$ has been chosen so that
\be
\lim_{x \to \infty}U_n(x;\alpha, \beta, N)=1. 
\ee
It is moreover manifest that on the real line $U_n(x;\alpha,\beta,N)$ has simple poles at $\alpha, \alpha+1, \dots, \alpha+n-1$ and hence we have:
\be
U_n(x;\alpha,\beta,N)= 1 + \sum_{k=0}^{n-1} \frac{\xi_{nk}}{x-\alpha -k} \lab{U_n_parfrac} \ee
for some coefficients $\xi_{nk}$. These rational functions of the ${_3}F_2$ type appear in the interpolation of the ratio of two Gamma functions $\frac{\Gamma(\alpha-\beta-z)}{\Gamma(\alpha-z)}$ by rational functions with prescribed poles and zeros (see \cite{Zhe_RIMS} for details).

The functions $U_n(x;\alpha,\beta,N)$ will be the center of attention in this paper. Using their explicit definition \re{HIF}, it is straightforward to verify that the operators $X^{(\alpha , \beta)}$, $Y^{(\alpha , \beta)}$ and $Z^{(\alpha , \beta)}$ all have the effect of 
producing the shift $\alpha \to \alpha+1$ in one of the parameters of $U_n(x;\alpha,\beta,N)$. More precisely we have:
\be
X^{(\alpha , \beta)} U_n(x;\alpha,\beta,N) = -\alpha \: U_n (x, \alpha+1,\beta,N), \lab{X_U} \ee
\be
Y^{(\alpha , \beta)} U_n(x;\alpha,\beta,N) = \alpha n(n+\beta-N) \: U_n (x, \alpha+1,\beta,N), \lab{Y_U} \ee
\be
Z^{(\alpha , \beta)} U_n(x;\alpha,\beta,N) = \frac{\alpha}{x-\alpha} \: U_n (x, \alpha+1,\beta,N). \lab{Z_U} \ee
This will prove an essential feature in the connection with  the bispectral properties of these rational functions of Hahn type. Before we look into that, we shall establish their biorthogonality.

\section{Biorthogonality}\lab{bior}
With hindsight we introduce the hypergeometric distribution:

\be
w^{(\alpha, \beta)}_x = \frac{(\beta -\alpha -N+2)_N}{(\beta -N+1)_N}\frac{(-N)_x (1-\alpha)_x}{x! (\beta-\alpha-N+2)_x}, \qquad \sum_{x=0}^N w_x=1\lab{w_x}, 
\ee
whose normalization is checked with the help of the Chu-Vandermonde formula: $_{2}F_1(-n,b;c;1)=(c-b)_n/(c)_n$. We can equip the space of rational functions $\mathcal{L}^{(\alpha, \beta)}$ involving two parameters $\alpha$ and $\beta$ and defined on the linear grid $x=0, 1, \dots, N$ with the scalar product:
\be
(f(x),g(x))_{(\alpha, \beta)}=\sum_{x=0}^N w^{(\alpha, \beta)}_x f(x)g(x),\qquad f, g \in  \mathcal{L}^{(\alpha, \beta)}.\lab{sp}
\ee
It proves useful to observe that the weight function $ w^{(\alpha, \beta)}_x$ is invariant under the following transformations denoted by $S$:
\be
S: \quad x \to (N-x),\quad \alpha \to (\beta +2-\alpha), \quad w^{(\alpha, \beta)}_x \to w^{(\beta +2 -\alpha, \beta)}_{(N-x)}=w^{(\alpha, \beta)}_x.\lab{S}
\ee
Now apply $S$ to $U_n(x;\alpha,\beta,N)$ to define the functions
\be
V_n(x;\alpha,\beta,N)=U_n(N-x;\beta +2-\alpha,\beta,N).
\ee
We note that the poles of these rational functions are located at $N+\alpha-\beta-2, N+\alpha-\beta-3, \dots, N+\alpha-\beta-n-1$.
We also perform the inversions $S$ on the operators $X$, $Y$ and $Z$ to obtain the transformed operators $\Tilde{X}$, $\Tilde{Y}$ and
$\Tilde{Z}$ respectively. Let $L$ be one of $X$, $Y$ or $Z$. Note that in view of the operation $x \to N-x$, in the case of operators, the symmetries $S$ must be accompanied by the exchange $T^+\leftrightarrow T^-$ in going from $L$ to $\Tilde{L}$.
One has the following identity:
\be 
\left(L^{(\alpha -1, \beta)}U_n(x;\alpha -1,\beta,N),V_m(x;\alpha,\beta,N)\right)=\left(U_n(x;\alpha,\beta,N),
\Tilde{L}^{(\alpha +1, \beta)}V_m(x;\alpha+1,\beta,N)\right),
\ee
which is readily proven by expressing the left hand side using the scalar product \re{sp} and then applying the transformations $S$ to recover the right hand side.

Given that the functions $U_n(x;\alpha,\beta,N)$ satisfy the GEVP 
$Y^{(\alpha , \beta)}U_n(x;\alpha,\beta,N)=
\lambda_nX^{(\alpha , \beta)}\\U_n(x;\alpha,\beta,N)$ 
with $\lambda_n=n(N-\beta-n)$, it is obvious that the functions $V_n(x;\alpha,\beta,N)$ will be solutions of
the GEVP $\Tilde{Y}^{(\alpha , \beta)}V_n(x;\alpha,\beta,N)=\lambda_n\Tilde{X}^{(\alpha , \beta)}V_n(x;\alpha,\beta,N)$ with the transformed operators and the same eigenvalues. A simple argument then shows that the sets of functions $\{U_n(x;\alpha,\beta,N)\}$ and $\{V_n(x;\alpha,\beta,N)\}$ are biorthogonal partners. Dropping the spectator parameters in $U_n$ and $V_n$, it goes like this:
\ba 
&&\left(Y^{(\alpha -1,\beta)}U_n(x;\alpha-1),V_m(x;\alpha)\right)
=\lambda_n\left(X^{(\alpha -1,\beta)}U_n(x;\alpha-1),V_m(x;\alpha)\right)\nonumber\\
&&=\left(U_n(x;\alpha),\Tilde{Y}^{(\alpha +1,\beta)}V_m(x;\alpha+1)\right)
=\lambda_m\left(U_n(x;\alpha),\Tilde{X}^{(\alpha +1,\beta)}V_m(x;\alpha+1)\right).
\ea
Observe as the notation suggests that the eigenvalues $\lambda_n$ do not depend on the parameter $\alpha$. It thus follows that
\be 
(\lambda_n -\lambda_m) \left((X^{(\alpha -1,\beta)}U_n(x;\alpha-1),V_m(x;\alpha)\right)=0.
\ee
Recalling from \re{X_U} that $X^{(\alpha -1,\beta)}U_n(x;\alpha-1)=-(\alpha -1)U_n(x;\alpha)$, this is seen to imply the biorthogonality relations:
\be 
\sum_{x=0}^N w^{(\alpha,\beta)}_x U_n(x;\alpha,\beta,N) V_m(x; \alpha, \beta,N) = 0, \quad \quad \mbox{if} \quad m \ne n \lab{biort} 
\ee
with $w^{(\alpha,\beta)}_x$ given by \re{w_x}. Thus, given the pair of operators $X^{(\alpha,\beta)}$ and $Y^{(\alpha,\beta)}$ we could construct explicitly the biorthogonal functions of Hahn type and obtained their orthogonality relations by considering the associated GEVP. We shall now proceed to show how GEVPs involving the triplet of operators $(X^{(\alpha,\beta)}, Y^{(\alpha,\beta)}, Z^{(\alpha,\beta)})$ characterize the bispectral properties of these functions. The tridiagonal action of these operators on the functions $U_n(x;\alpha,\beta,N)$ which we present next will be key.

\section{Tridiagonal actions of the triplet $X,Y,Z$}
We already know that the operators $X^{(\alpha,\beta)},Y^{(\alpha,\beta)},Z^{(\alpha,\beta)}$ are tridiagonal matrices \re{X_e} -\re{Y_e} with respect to the basis $e_n(x)$  defined on the finite uniform grid $x=0,1,\dots,N$. More exactly, the operator $Y{(\alpha,\beta)}$ is tridiagonal, while the operators $X{(\alpha,\beta)}$ and $Z{(\alpha,\beta)}$ are bidiagonal (the upper off-diagonal is zero).

In turn, in the basis $\varphi_n(x)$ these matrices are tridiagonal too and are given by \re{XM_phi} - \re{ZM_phi} (again only the operator $Y$ is tridiagonal, while $X$ and $Z$ are bidiagonal).

Consider anew the action of these operators on the rational functions $U_n(x;\alpha,\beta,N)$. The functions $U_n(x)$ form another basis on the finite-dimensional space $M_N$. From the contiguity relations of the ${_3}F_2(1)$ series in particular, one finds that \re{X_U}, \re{Y_U} and \re{Z_U} can be reexpressed as follows:
\begin{align}
& X^{(\alpha,\beta)} U_n(x;\alpha,\beta,N) = \mu^{(1)}_n U_{n+1}(x;\alpha,\beta,N) + \mu_n^{(2)} U_n(x;\alpha,\beta,N) + \mu_n^{(3)} U_{n-1}(x;\alpha,\beta,N), \lab{XU} \\
& Y^{(\alpha,\beta)} U_n(x;\alpha,\beta,N) = \mu^{(4)}_n U_{n+1}(x;\alpha,\beta,N) + \mu_n^{(5)} U_n(x;\alpha,\beta,N) + \mu_n^{(6)} U_{n-1}(x;\alpha,\beta,N) \lab{YU}
\end{align}
and
\be
Z^{(\alpha,\beta)} U_n(x;\alpha,\beta,N) = \mu^{(7)}_n U_{n+1}(x;\alpha,\beta,N) + \mu_n^{(8)} U_n(x;\alpha,\beta,N) + \mu_n^{(9)} U_{n-1}(x;\alpha,\beta,N), \lab{ZU} \ee
where
\ba
&\displaystyle
 \mu_n^{(1)} = {\frac {n \left( \beta+n+1 \right)  \left( N-\beta-n \right) }{
 \left( N-\beta-2\,n \right)  \left( N-\beta-2\,n-1 \right) }}, \lab{mu1} 
\\
&\displaystyle
\mu_n^{(2)} = -\alpha-n+{\frac {\beta\,n \left( n-1 \right) +{n}^{2} \left( n-1
 \right) }{N-\beta-2\,n+1}}-{\frac {\beta\,n \left( n+1 \right) +n
 \left( n+1 \right) ^{2}}{N-\beta-2\,n-1}}, \lab{mu2} 
\\
&\displaystyle
 \mu_n^{(3)} = {\frac {n \left( N-\beta-n \right)  \left( N-n+1 \right) }{ \left( N-
\beta-2\,n \right)  \left( N-\beta-2\,n+1 \right) }}, \lab{mu3}
\\
&\displaystyle
 \mu_n^{(7)} = -{\frac { \left( \beta+n+1 \right)  \left( N-\beta-n \right) }{
 \left( N-\beta-2\,n \right)  \left( N-\beta-2\,n-1 \right) }}, \lab{mu7} 
\\
&\displaystyle
 \mu_n^{(8)}={\frac { \left( n+1 \right) \beta+ \left( n+1 \right) ^{2}}{N-\beta-2
\,n-1}}-{\frac {\beta\,n+{n}^{2}}{N-\beta-2\,n+1}},
  \lab{mu8} 
\\
&\displaystyle
 \mu_n^{(9)}=-{\frac {n \left( N-n+1 \right) }{ \left( N-\beta-2\,n \right) 
 \left( N-\beta-2\,n+1 \right) }.} \lab{mu9} 
\ea


As for the coefficients $\mu^{(4)}_n, \: \mu^{(5)}_n, \: \mu^{(6)}_n$, they are proportional to the coefficients $\mu^{(1)}_n, \: \mu^{(2)}_n, \: \mu^{(3)}_n$ as per equations \re{X_U}-\re{Y_U}:
\be
\mu^{(3+i)}_n = n(N-n-\beta) \: \mu_n^{(i)}, \quad i=1,2,3 .\lab{muY} 
\ee 
In the formulas above defining the parameters $\mu^{(\ell)}_n, \ell = 1,...,8$, the range of $n$ is $\{0, \dots ,N\}$ except for $\ell = 1, 4, 7$ for which cases the value $N$ is excluded and the expressions replaced by 
\be
\mu^{(1)}_N = \mu^{(4)}_N = \mu^{(7)}_N =0. \lab{mu=0} \ee
Condition \re{mu=0} provides the finite-dimensional restrictions of the operators $X,Y,Z$ in the basis $U_n(x)$

Equivalently, in the basis $U_n(x)$ the operators $X,Y,Z$ can be presented by square matrices of size $N+1$. For example
\be
X=\left[
\begin{array}{ccccc}
\mu_0^{(2)} & \mu_1^{(3)} & 0 & \cdots & 0 \\
\mu_0^{(1)} & \mu_1^{(2)} & \mu_2^{(3)} &  & \vdots  \\
0 & \ddots & \ddots & \ddots & 0 \\
\vdots & & \mu_{N-2}^{(1)} & \mu_{N-1}^{(2)} & \mu_N^{(3)} \\
0 & \cdots & 0 & \mu_{N-1}^{(1)} & \mu_N^{(2)}
\end{array}
\right],
\lab{XU_phi} \ee
and similar expressions for the matrices $Y$ and $Z$.

\section{Bispectrality}
We just saw that the operators $X^{(\alpha,\beta)},Y^{(\alpha,\beta)},Z^{(\alpha,\beta)}$ are all tridiagonal in the basis $U_n(x;\alpha,\beta,N)$ and hence can be represented by corresponding tridiagonal square matrices of size $N+1$. This observation leads directly to the bispectrality of the rational functions $U_n(x;\alpha,\beta,N)$. This means that these functions satisfy a pair of GEVPs, one in the variable $x$ with the eigenvalue depending on $n$, referred to as the difference equation, the other in the degree $n$ with the eigenvalue depending on $x$, referred to as the recurrence relation.

\subsection{Difference equation}

Consider again the GEVP for the pair of the operators $X^{(\alpha,\beta)}$ and $Y^{(\alpha,\beta)}$. We already know from Section \re{secrat} that the functions $U_n(x;\alpha,\beta,N)$ are solutions of 
\be 
Y^{(\alpha , \beta)}U_n(x;\alpha,\beta,N)=\lambda_nX^{(\alpha , \beta)}\\U_n(x;\alpha,\beta,N)\lab{GEVP_XY}
\ee
with the eigenvalues given by
\be
\lambda_n = n(N-n-\beta) \lab{lambda_n}. 
\ee
Note that this fact is also borne out by the equations \re{X_U} and \re{Y_U}.

Using the explicit expressions \re{X_op} and \re{Y_op} of the operators $X^{(\alpha,\beta)}$ and $Y^{(\alpha,\beta)}$, the equation \re{GEVP_XY} reads:
\ba
&& A_1(x) U_n(x+1;\alpha,\beta,N) + A_2(x) U_n(x-1;\alpha,\beta,N) + A_3(x) U_n(x;\alpha,\beta,N) = \nonumber \\
&& \lambda_n \left[ (x-\alpha) U_n(x;\alpha,\beta,N) -x U_n(x-1,;\alpha,\beta,N) \right]. \lab{GEVP_x} 
\ea
This yields a second-order difference equation for the rational functions $U_n(x;\alpha,\beta,N)$.

\subsection{Recurrence relation}
Take now the pair of the operators $X^{(\alpha , \beta)}$ and $Z^{(\alpha , \beta)}$. From \re{X_U} and \re{Y_U} it follows that the biorthogonal rational functions $U_n(x;\alpha,\beta,N)$ satisfy another GEVP namely,

\be
X^{(\alpha , \beta)} U_n(x) = (\alpha-x) Z^{(\alpha , \beta)}U_n(x;\alpha,\beta,N). \lab{GEVP_XZ} 
\ee
Using \re{XU} and \re{ZU}, this equation becomes
\ba
&& \mu^{(1)}_n U_{n+1}(x;\alpha,\beta,N) + \mu_n^{(2)} U_n(x;\alpha,\beta,N) + \mu_n^{(3)} U_{n-1}(x;\alpha,\beta,N) = \nonumber \\
&& (\alpha-x) \left[ \mu^{(7)}_n U_{n+1}(x;\alpha,\beta,N) + \mu_n^{(8)} U_n(x;\alpha,\beta,N) + \mu_n^{(9)} U_{n-1}(x;\alpha,\beta,N) \right] \lab{GEVP_n} \ea
This is recognized as a 3-term recurrence relation.

We thus see that both the recurrence relation and the difference equation of the rational functions $U_n(x;\alpha,\beta,N)$ arise as GEVPs
associated to the two pairs of operators: $(X^{(\alpha , \beta)},Y^{(\alpha , \beta)})$ and $(X^{(\alpha , \beta)},Z^{(\alpha , \beta)})$. This
is how the triplet $\{X^{(\alpha , \beta)}, Y^{(\alpha , \beta)}, Z^{(\alpha , \beta)}\}$ accounts for the bispectrality of the biorthogonal rational functions of Hahn type.
\subsection{Comparison with classical orthogonal polynomials}
It is common to call classical, functions that are solutions of a bispectral problem. Let us stress the striking differences in the expression of bispectrality that we are observing between the biorthogonal rational functions of Hahn type and the classical orthogonal polynomials (OPs). In the case of OPs $\{P_n(x)\}$, there exists a pair of operators $X$ and $Y$ bearing the name of Leonard (in the finite-dimensional case), such that $X$ is diagonal with respect to the basis $e_k(s)$ and that $Y$ is tridiagonal in that same basis (see \cite{Ter} for details). This provides the second order difference equation 
\be
Y P_n(x(s)) = A_1(s) P_n(x(s+1)) + A_2(s) P_n(x(s)) + A_3(s) P_{n-1}(x(s-1)) = \lambda_n P_n(x(s)) \lab{diff_P} 
\ee
where $x(s)$ is the orthogonality grid.  
In the basis $P_n(x)$, in contrast, the operator $X$ is tridiagonal and $Y$ is diagonal. This gives the recurrence relation that ensures orthogonality:
\be
X P_n(x) = P_{n+1}(x) + b_n P_n(x) + u_n P_{n-1}(x). \lab{P_rec} \ee 
For the rational biorthogonal functions we have considered, there is a triplet of operators $(X,Y,Z)$ instead of a Leonard pair. Moreover, it is seen that in both the bases $e_n(s)$ and $U_n(x;\alpha,\beta,N)$ each one of the operators $X,Y,Z$ is tridiagonal. Importantly, in the basis $e_n(x)$, the difference operators $X^{(\alpha , \beta)}$ and $Z^{(\alpha , \beta)}$ are related to one another by the factor $\alpha-x$ which depends only on $x$; while in basis $U_n(x;\alpha,\beta,N)$ the operators $X^{(\alpha , \beta)}$ and $Y^{(\alpha , \beta)}$ are connected by the factor  $\lambda_n=n(N-n-\beta)$ which depends only on $n$. This feature of the triplet $\{X,Y,Z\}$ associated to the classical rational functions of Hahn type can be viewed as an analog of the property (see above) of the Leonard pairs attached to classical orthogonal polynomials. One will bear in mind that the bispectral rational functions obey a couple of GEVPs while the classical OPs satisfy two regular eigenvalue equations.

Finally, let us underscore that the basis $\varphi_n(x;\alpha)$ introduced in Section \re{secrat} plays for the rational functions of Hahn type a role comparable to the split basis in the OP context \cite{Ter}. In such bases the operators $X$ and $Y$ of a Leonard pair are bidiagonal and this allows to obtain the explicit hypergeometric expressions of the corresponding orthogonal polynomials.
We observed that the operator $Y^{(\alpha , \beta)}$ given in \re{Y_op} remains tridiagonal in the basis $\varphi_n(x;\alpha)$, however due to the factorization property \re{Y-Y1}, solving the GEVP  \re{GEVP_XY} again amounts to obtaining the solution of a two-term recurrence relation.

Summing up we have found that the operators $X^{(\alpha,\beta)}, Y^{(\alpha,\beta)}, Z^{(\alpha,\beta)}$ are tridiagonal (or bidiagonal) in three  different bases: $\{e_n(x)\}$, $U_n(x;\alpha,\beta,N)$ and $\{\varphi_n(x;\alpha)\}$. Owing to special factorization properties, the GEVP \re{GEVP_XY} and  \re{GEVP_XZ} either yield the explicit hypergeometric expression of these biorthogonal rational functions or provide their difference equation and recurrence relation.

\section{The quadratic algebra realized by $X^{(\alpha,\beta)}, Y^{(\alpha,\beta)}, Z^{(\alpha,\beta)}$}
The algebras realized by the recurrence and difference operators of the orthogonal hypergeometric operators are said to be of the Askey-Wilson type \cite{GLZ_Annals} and have proven ubiquitous and of prime interest. They encode the bispectral properties of the corresponding orthogonal polynomials which are determined by the representation theory of these algebras. It is to be expected that the relations obeyed by the operators
$X^{(\alpha,\beta)}, Y^{(\alpha,\beta)}, Z^{(\alpha,\beta)}$ define an algebra which we will denote $\mathcal{R}_H$ that will similarly embody most features of the biorthogonal rational functions of Hahn type. We provide these relations below.  For conciseness, we omit specifying the parameters. The relations are obtained by direct computation using the definitions \re{X_op}, \re{Z_op}, \re{Y_op} and read:
\begin{eqnarray}
\mbox{}&[Z,X] = Z^2 +Z, \lab{CMZX} \\%
\mbox{}&[X,Y] = \xi_1(X^2+Z^2) + \{X,Z\} +\{Y,Z\} + \xi_2 X + \xi_3 Z +Y + \xi_0 \mathcal{I}, \lab{CMXY} \\%
\mbox{}&[Y,Z] = 3 X^2 + Z^2 +\xi_1 \{X,Z\} +\xi_4 X + \xi_2 Z + \xi_0 \mathcal{I}, \lab{CMYZ}
\end{eqnarray}
where $\{X,Y\} \equiv XY+YX$ stands for the anticommutator, $\mathcal{I}$ is the identity operator and where
\be
\xi_1= N-\beta, \: \xi_0 = \alpha(\alpha-\beta), \; \xi_2= 1+\alpha \left(1+N-\beta  \right), \; \xi_3=\alpha+(\alpha+1)(N-\beta), \: \xi_4 =4 \alpha-2 \beta-1. \lab{xi} \ee
Manifestly this is quadratic algebra since all commutators are expressed as second order polynomials in the non-commutative variables $X,Y,Z$.  

It can be verified that the Casimir operator of the algebra thus defined without the restrictions \eqref{xi} is the following cubic polynomial in the generators $X,Y,Z$:
\begin{align}
Q=& X^3 + \xi_1 Z^3 + \frac{\xi_1}{2} \{X^2,Z\} + 2\: \{X,Z^2\} + \frac{1}{2} \: \{Y, Z^2\} + \nonumber\\
& \frac{\xi_4}{2} X^2 + \frac{\xi_1+\xi_3+\xi_4}{2}\:  Z^2 + \frac{\xi_2+3}{2} \: \{X,Z\}  +   \frac{1}{2}\: \{Y,Z\} + (\xi_0+1/2)\: X +    (\xi_0 + \xi_4/2) \: Z. \lab{Q}
\end{align}
In the realization provided by the difference operators \re{XU}-\re{ZU}, this Casimir operator is proportional to the identity operator and becomes precisely:
\be
Q= \frac{\alpha(\beta-\alpha)}{2} \mathcal{I}. \lab{Q_val} 
\ee

\remark We noted in Section \re{bior} that the operators $\Tilde{X}^{(\alpha,\beta)}$, $\Tilde{Y}^{(\alpha,\beta)}$ and $\Tilde{Z}^{(\alpha,\beta)}$ which account through GEVPs for the bispectrality of the biorthogonal partners $V_n(x;\alpha,\beta,N)$ of the rational functions $U_n(x;\alpha,\beta,N)$, are obtained from the operators $X^{(\alpha,\beta)}$, $Y^{(\alpha,\beta)}$ and $Z^{(\alpha,\beta)}$ under the symmetry operation $S$ given in \re{S} supplemented by the exchange $T^+ \leftrightarrow T^-$. It follows as expected, that the operators $\t X, \t Y, \t Z$ satisfy the relations \re{CMZX}-\re{CMYZ} of the same quadratic algebra $\mathcal{R}_H$ with $\alpha$ replaced by $\alpha-\beta-2$. in the parameters $\xi_i$.

\remark We will also point out that the algebra presented above can be derived from a potential and share in that a property of most Calabi-Yau algebras \cite{G}, \cite{BRS} of dimension 3. Let $F=\mathbb{C}[x_1, x_2, \dots x_n]$ be a free associative algebra with $n$ generators and $F_{cycl}=F/[F,F]$. $F_{cycl}$ has the cyclic words $[x_{i_1}, x_{i_2}\dots x_{i_r}]$ as basis; these can be viewed as oriented closed paths. The map $\frac{\partial}{\partial x_j}: F_{cycl} \to F$ is such that
\be 
\frac{\partial [x_{i_1}, x_{i_2},\dots x_{i_r}] }{\partial x_j}= \sum_{\{s|i_s=j\}} x_{i_s +1}x_{i_s +2}\dots x_{i_r}x_{i_1}x_{i_2}\dots x_{i_s -1} 
\ee
and extended to  $F_{cycl}$ by linearity on combinations of cyclic words. Let $\Phi(x_1,\dots x_n) \in F_{cycl}$. An algebra whose defining relations are given by 
\be 
\frac{\partial \Phi}{\partial x_j}=0, \qquad j=1,\dots n
\ee
is said to derive from the potential $\Phi$. Now let $x_1=X$, $x_2=Y$ and $x_3=Z$ and take
\ba 
&&\Phi=[XYZ]-[YXZ]-[X^3]-\xi _1[X^2Z]-[XZ^2]-[YZ^2] \nonumber\\
&&\phantom{\Phi=}-\frac{\xi_1}{3}[Z^3]-\frac{\xi _4}{2}[X^2
]-\xi_2[XZ]-[YZ]-\frac{\xi _3}{2}[Z^2]-\xi _0([X]+[Z]).\lab{Pot}
\ea
With this, it is not hard to see that the relations \re{CMZX}, \re{CMXY} and \re{CMYZ} of $\mathcal{R}_H$ are respectively given by 
$\frac{\partial \Phi}{\partial Y}=0$, $\frac{\partial \Phi}{\partial Z}=0$, $\frac{\partial \Phi}{\partial X}=0$ and that the algebra associated to the bispectral properties of the rational functions of Hahn type derives from the potential \re{Pot}.

\section{Conclusion}

In summary, we have found that the bispectral properties of the rational biorthogonal functions of Hahn type are described by a triplet of operators $X, Y, Z$. These operators are tridiagonal in three bases naturally associated to the context: the canonical basis on a linear grid, the ratios of Pochhammer symbols and the rational functions themselves. These operators all induce the same shift in one parameter of the rational functions accompanied by factors that are either constant, involving the variable only or depending solely on the degree. This entails the bispectrality of the rational functions. Indeed,  as a result of these factorizations and of the tridiagonal actions, the difference equation and the recurrence relation emerge simply from generalized eigenvalue problems associated with two pairs of operators out of the triplet $X, Y, Z$. This offers a parallel in the case of biorthogonal rational functions to the description of bispectrality in terms of Leonard pairs for orthogonal polynomials. Of note in addition is the fact that factorization also occurs in the basis consisting of ratios of shifted factorials so that the associated GEVP can be solved directly and the rational functions explicitly obtained much like the bidiagonal action of Leonard operators in the split basis permits the determination of the associated orthogonal polynomials from a two-term recurrence relation. 

It would be very pertinent to construct the representations of the 
algebra $\mathcal{R}_H$ that the triplet $X, Y, Z$ realizes; this is a project we have in mind as these representations should encompass the characterization of the rational functions of Hahn type.

We trust that this general framework extends, beyond the case we have studied, to other bispectral rational functions. It would hence be relevant to explore how it is realized for other biorthogonal rational functions of interest such as those of Wilson \cite{W} or their specializations \cite{GM}. (A central question is the identification of tridiagonal operators that admit factorization instead of the Leonard mutual diagonalization.) We certainly plan on doing that with an eye to identifying the associated algebras and elaborating a general representation theoretic description of biorthogonal rational functions that are bispectral.

\bigskip\bigskip
{\Large\bf Acknowledgments}
The authors are grateful to G. Bergeron and J. Gaboriaud for discussions and in particular for verifying abstractly the expression of the Casimir element.
S.T. wishes to thank the Centre de Recherches Math\'ematiques for its hospitality during the course of this investigation.
His work is partially supported by JSPS KAKENHI (Grant 
Numbers 19H01792, 17K18725). The research of L.V. is funded in part by a
Discovery Grant from the Natural Sciences and Engineering Council (NSERC) of Canada. 
A. Z. gratefully acknowledges the award of a CRM-Simons Professorship and
is supported by the National Science
Foundation of China (Grant No.11771015).

\vspace{5mm}

\bb{99}
\bi{W}J. A. Wilson, {\it Orthogonal functions from Gram determinants}, SIAM J. Math. Anal. {\bf 22} (1991), 1147-1155

\bi{GM}D. P. Gupta and D. R. Masson, {\it Contiguous relations, continued fractions and orthogonality}, Trans. Amer.
Math. Soc. {\bf 350} (1998), 769-808.

\bi{IM} M. E. H. Ismail and D. R. Masson, {\it Generalized orthogonality and continued fractions}, J. Approx. Theory
{\bf 83} (1995), 1-40.

\bi{Zhe_RIMS} A. Zhedanov, {\it Pad\'e interpolation table
and biorthogonal rational functions}, Rokko. Lect. in Math. {\bf 18} (2005), 323-363.

\bi{Zhe_GEVP} A. Zhedanov, {\it Biorthogonal rational functions and the generalized eigenvalue problem},
J. Approx. Theory {\bf 101} (1999), 303-329.

\bi{KM} E. G. Kalnins and W. Miller Jr, {\it q-series and orthogonal polynomials associated with Barnes' first lemma}, SIAM J. Math. Anal. {\bf19} (1988), 1216-1231

\bibitem{KLS} R. Koekoek, P. A. Lesky, and R. F. Swarttouw. {\it Hypergeometric orthogonal polynomials and their q-analogues}. Springer, 1-st edition, 2010.

\bi{BI} E. Bannai and T. Ito
  {\it Algebraic combinatorics},
Benjamin/Cummings,1984.

\bi{GLZ_Annals} Ya. A. Granovskii, I. M. Lutzenko, and A. Zhedanov, {\it Mutual integrability, quadratic algebras,
and dynamical symmetry}. Ann. Phys. {\bf 217} (1992),  1-20.

\bi{Ter} P. Terwilliger, {\it Two linear transformations each tridiagonal with respect to an eigenbasis of the other}, Lin. Alg. Appl. {\bf 330} (2001), 149-203.

\bi{L} D. A. Leonard, {\it Orthogonal polynomials, duality and association schemes}, SIAM J. Math. Anal. {\bf 13} (1982), 656-663.

\bi{TVZ} S. Tsujimoto, L. Vinet and A. Zhedanov, {\it The rational Heun operator and Wilson biorthogonal functions}, arXiv preprint (2019) arXiv:1912.11571.

\bi{VZ} L. Vinet and A. Zhedanov, {\it The Heun operator of Hahn type}, 
Proc. Amer. Math. Soc.
{\bf 147} (2019),
2987-2998

\bi{G} V. Ginzburg, {\it Calabi-Yau algebras}, arXiv preprint (2006), arXiv: math/0612139

\bi{BRS} G. Bellamy, D. Rogalski, T. Schedler, J. T. Stafford, and M. Wemyss, {\it Noncommutative algebraic geometry}, Cambridge University Press, 2016.
\end{thebibliography}

\end{document}